\documentclass[a4paper,reqno]{amsart}
\usepackage{amsfonts}
\theoremstyle{plain}
\begingroup
\newtheorem{theorem}{Theorem}[section]
\newtheorem{lemma}[theorem]{Lemma}
\newtheorem{proposition}[theorem]{Proposition}

\endgroup

\theoremstyle{definition}
\begingroup
\newtheorem{definition}[theorem]{Definition}

\newtheorem{example}[theorem]{Example}
\endgroup

\theoremstyle{remark}
\begingroup
\endgroup

\numberwithin{equation}{section}

\mathchardef\emptyset="001F

\newcommand{\Om}{\Omega}
\newcommand{\R}{{\mathbb R}}
\newcommand{\Rtwo}{{\R}^2}
\newcommand{\wto}{\rightharpoonup}
\newcommand{\setmeno}{\!\setminus\!}
\newcommand{\hone}{{\mathcal H}^1}
\newcommand{\F}{{\mathcal F}}
\newcommand{\G}{{\mathcal G}}
\newcommand{\subsethn}{\mathrel{\mathop{\smash\subset
\vphantom{=}}\limits^\sim} }
\newcommand{\eqhn}{\mathrel{\mathop{\smash=
\vphantom{\scriptscriptstyle a}}\limits^\sim}}

\title[Unilateral slope for the Mumford-Shah functional]
{On a notion of unilateral slope \\ for the Mumford-Shah functional}
\author[Gianni Dal Maso]{Gianni Dal Maso}
\address[Gianni Dal Maso]{SISSA, Via Beirut 2, 34014 Trieste,
Italy}
\email[Gianni Dal Maso]{dalmaso@sissa.it}
\author[Rodica Toader]{Rodica Toader}
\address[Rodica Toader]{Dipartimento di Ingegneria Civile, Via delle
Scienze 208, 33100
  Udine, Italy}
\email[Rodica Toader]{toader@uniud.it}

\begin{document}
\begin{abstract}
In this paper we introduce a notion of unilateral slope for the
Mumford-Shah functional, and provide an explicit formula in the case of
smooth cracks. We show
that the slope is not lower semicontinuous and study the corresponding
relaxed functional.
\end{abstract}
\maketitle
{\small

\bigskip

\bigskip
}
\bigskip
\bigskip
\begin{section}{INTRODUCTION}

In the study of the gradient flow for a nonsmooth functional $\F$ on a
metric space $(X,d)$, it is useful to introduce a notion of slope
$|\partial\F|$, which coincides with the norm of the gradient
$\nabla\F$  in the case of a smooth functional on a Hilbert space.
For every $u\in X$ with $\F(u)<+\infty$ the slope $|\partial\F|(u)$ is
defined  by
$$
|\partial\F|(u):=\limsup_{v\to u}\frac{(\F(u)-\F(v))^+}{d(v,u)}\,,
$$
where $(\cdot)^+$ denotes the positive part. For the general properties
of the slope and for the comparison with other classical notions we
refer to \cite{DG-Mar-Tos}, \cite{Deg-Mar-Tos}, \cite{Mar-Sac-Tos},
\cite{Amb-Mov}, \cite{Amb-Gig-Sav} and to the forthcoming book
\cite{Amb-Gig-Sav-b}.

In this paper we begin the study of this notion for the Mumford-Shah
functional defined in the space  $SBV(\Om)$ of special functions with
bounded variation (see~\cite{A-F-P}) on a bounded open set $\Om\subset\Rtwo$ with $C^1$ boundary.  

In view of the applications to irreversible crack growth in
fracture mechanics, it is convenient to write this functional in the
form
\begin{equation}\label{MS}
\F(u,S):=\|\nabla u\|^2 + \hone(S)\,,
\end{equation}
using two independent variables: the function  $u\in SBV(\Om)$, which plays the
role of the displacement, and the set $S$, which plays the role of the crack. Here and henceforth $\|\cdot\|$ denotes the $L^2$ norm,  $\nabla$ denotes the (approximate) gradient, and $\hone$ is the one dimensional Hausdorff measure. Since the displacement must be (approximately) continuous out of the
crack, the domain of the functional $\F(u,S)$ satisfies the constraint
$S(u)\subsethn S$, where $S(u)$ is the jump set of $u$ and $\subsethn$
denotes inclusion up to an $\hone$-negligible set. For the precise definitions of all these notions we refer to \cite{A-F-P}. 

The irreversibility of crack growth leads to the following unilateral
variant of the notion of slope:
\begin{equation}\label{slope}
|\partial\F|(u,S):=\limsup_{v\to u}\frac{(\F(u,S)-\F(v,S\cup
S(v)))^+}{\|v-u\|}\,,
\end{equation}
where $v\to u$ in $L^2(\Om)$. We consider also the case with fixed
boundary conditions on $\partial\Om$
\begin{equation}\label{slopeb}
|\partial\F|_{b}(u,S):=\limsup_{{\scriptstyle v\to  u\atop  \scriptstyle v=u \text{ on }\partial\Om}}\frac{(\F(u,S)-\F(v,S\cup
S(v)))^+}{\|v-u\|}\,,
\end{equation}
where the equality $v=u$ on $\partial\Om$ means that the traces of $v$ and $u$ agree $\hone$-a.e.\ on $\partial\Om$. It is easy to see that $|\partial\F|_{b}(u,S)\le |\partial\F|(u,S)$.

For every $S\subset\Om$ with $\hone(S)<+\infty$ we define $SBV^2(\Om,S)$
as the set of functions $u\in SBV(\Om)\cap L^2(\Om)$ with $\nabla u\in
L^2(\Om;\R^2)$ and $S(u)\subsethn S$; we define $SBV^2_0(\Om,S)$ as the
set of functions $v\in SBV^2(\Om,S)$ whose trace vanishes $\hone$-a.e.\
on $\partial\Om$.
We use the symbol $(\cdot|\cdot)$ to denote the scalar product in $L^2(\Om)$ or
$L^2(\Om;\R^2)$, according to the context.

A necessary condition for the finiteness of the slope is given by the
following proposition. 

\begin{proposition}\label{necessary}
Let $S$ be a subset of $\Om$ with $\hone(S)<+\infty$ and let $u\in
SBV^2(\Om,S)$.
If $|\partial\F|_b(u,S)<+\infty$, then there exists $f\in L^2(\Om)$
such that
\begin{equation}\label{nec2}
(\nabla u| \nabla\varphi)=(f|\varphi)\qquad\forall \varphi\in
SBV^2_0(\Om,S)\,.
\end{equation}
Moreover, $2\|f\|\leq |\partial\F|_b(u,S)$.

If, in addition, $|\partial\F|(u,S)<+\infty$, then we have also
\begin{equation}\label{nec1}
(\nabla u| \nabla\varphi)=(f|\varphi)\qquad\forall \varphi\in
SBV^2(\Om,S)\,.
\end{equation}
\end{proposition}

If $S$ and $u$ are sufficiently smooth, condition (\ref{nec2}) is
equivalent to say that $-\Delta u=f$ in $\Om\setmeno S$ and $\partial
u/\partial n=0$ on $S$, while (\ref{nec1}) implies also that $\partial
u/\partial n=0$ on $\partial\Om$.

For every open set $U\subset\Rtwo$ the space of distributions on $U$ is
denoted by  ${\mathcal D}'(U)$.
The following theorem shows that the slope in the previous proposition
is given exactly by $2\|f\|$ when $S$ and $u$ are sufficiently smooth,
and $u$ satisfies the Neumann boundary condition.

\begin{theorem}\label{smooth}
Let $S$ be a one dimensional $C^1$ manifold without boundary contained
in $\Om$ and let $u$ be a function such $u|_{\Om_i}\in
C^1(\overline\Om_i)$ for every connected component $\Om_i$ of
$\Om\setmeno S$. Assume that $-\Delta u=f$ in ${\mathcal D}'(\Om\setmeno S)$,
with $f\in L^2(\Om)$, and that $\partial u/\partial n=0$ on $S$. Then
$|\partial\F|_{b}(u,S)=2\|f\|$. If, in addition, $\partial u/\partial
n=0$ on $\partial\Om$, then $|\partial\F|(u,S)=2\|f\|$.
\end{theorem}

The hypothesis that $S$ is a manifold without boundary is crucial in the previous theorem. In Example~\ref{example1} we will show that, if $S$ is a closed segment in $\Om$ and $u$ is harmonic on $\Om\setmeno S$, satisfies the Neumann boundary condition on $S$, and has a large stress intensity factor at one of the crack tips, then $|\partial\F|_b(u,S)>0$. On the other hand we will show in Example~\ref{example2} that, in this case, there exists a sequence $S_k$ of $C^1$ manifolds without boundary, which converges to $S$ in any reasonable sense, and a sequence $u_k$ of harmonic functions on $\Om\setmeno S_k$, satisfying the Neumann boundary condition on $S_k$, and such that $u_k$ converges to $u$ strongly in $L^2(\Om)$ and $\nabla u_k$ converges to $\nabla u$  strongly in $L^2(\Om;\R^2)$. By Theorem~\ref{smooth} this implies that 
$|{\partial\F}|_b(u_k,S_k)=0$ and shows that $|{\partial\F}|_b$ is not lower semicontinuous.

In the general theory of gradient flows, when the slope is not lower semicontinuous, its lower semicontinuous envelope plays an important role. In the
case of the Mumford-Shah functional, for every $S\subset\Om$ with
$\hone(S)<+\infty$ and every $u\in SBV^2(\Om,S)$,  the natural notion of
lower semicontinuous envelope is defined by
$$
|\overline{\partial\F}|(u,S):=\inf\{\liminf_{k\to\infty}
|\partial\F|(u_k,S_k)\}\,,
$$
where the infimum is taken over all sequences $(u_k,S_k)$ such that
$u_k\to u$ strongly in $L^2(\Om)$, $\nabla u_k\wto \nabla u$ weakly in
$L^2(\Om;\Rtwo)$, and $S_k$ $\sigma^2$-converges to $S$, according to
Definition~\ref{sigma2} below.
Similarly we define
$$
|\overline{\partial\F}|_b(u,S):=\inf\{\liminf_{k\to\infty}
|\partial\F|_b(u_k,S_k)\}\,,
$$
where the infimum is taken over the same set of sequences.

We prove the following general result on the relaxed slope.

\begin{proposition}\label{relaxed}
Let $S$ be a subset of $\Om$ with $\hone(S)<+\infty$ and let $u\in
SBV^2(\Om,S)$. If  $|\overline{\partial\F}|_b(u,S)<+\infty$, then there
exists $f\in L^2(\Om)$ such that
\begin{eqnarray}
& -{\rm div}(\nabla  u)=f\quad\hbox{ in }{\mathcal D}'(\Om)\,,
\label{div}
\\
&
|\nabla u|^2-{\rm div}(u\nabla  u)\leq fu \quad\hbox{ in }{\mathcal
D}'(\Om)\,.
\label{divu}
\end{eqnarray}
Moreover,
$2\|f\|\leq |\overline{\partial\F}|_b(u,S)$.

If, in addition,  $|\overline{\partial\F}|(u,S)<+\infty$, then we have also
\begin{equation}\label{divtilde}
-{\rm div}(\nabla  \tilde u)=\tilde f\quad\hbox{ in }{\mathcal
D}'(\Rtwo)\,,
\end{equation}
where the tilde denotes the zero extension to $\Rtwo$.
\end{proposition}

If $S$ and $u$ are sufficiently smooth, condition (\ref{div}) is
equivalent to say that $-\Delta u=f$ in $\Om\setmeno S$ and $\partial
u/\partial n$ is continuous across $S$, while (\ref{divtilde}) implies
also that $\partial u/\partial n=0$ on $\partial\Om$. Note that
(\ref{nec2}) implies the stronger condition $\partial u/\partial n=0$
on $S$. Condition (\ref{divu}) says that $ [u]\,\partial u/\partial
n\ge 0$ on $S$, where $[u]$ denotes the jump of $u$ on $S$.

The following theorem shows that in the previous proposition we have
$|\overline{\partial\F}|_b(u,S)=2\|f\|$ if $S$ is partially smooth,
and $u$ satisfies (\ref{div}) and (\ref{divu}) with an $f\in L^p(\Om)$
with $p>2$.

\begin{theorem}\label{thmrelax}
Assume that $\Om$ has a $C^2$ boundary and 
let $S$ be a compact subset of $\Om$. Suppose that
there exists a finite set
$F\subset S$ such that $S\setmeno F$  is a one dimensional $C^1$ manifold, and that for every $x\in F$ there exists an open neighbourhood $U$ of $x$ such that $U\setmeno S$ can be expressed as the union of a finite number of open sets with Lipschitz boundary and of a set of Lebesgue measure zero.
Let
$u\in SBV^2(\Om,S)\cap L^\infty(\Om)$ be a function which satisfies (\ref{div}) and (\ref{divu}) with
$f\in L^p(\Om)$, $p>2$. Assume that $[u]\neq0$ $\hone$ a.e.\ on $S$. Then
$|\overline{\partial\F}|_b(u,S)=2\|f\|$.
If, in addition, (\ref{divtilde}) holds, then
$|\overline{\partial\F}|(u,S)=2\|f\|$.
\end{theorem}

The main difference between Theorems~\ref{smooth} and~\ref{thmrelax} is
that in Theorem~\ref{smooth} the set $S$ is assumed to be  smooth,
while in
Theorem~\ref{thmrelax} it is smooth except possibly for a finite number
of points. The comparison between these results shows that the slope is sensitive to the
behaviour of $u$ near the crack tips, while this is not the case for the
relaxed slope. Another difference is the fact that in
Theorem~\ref{smooth} we assume the Neumann condition $\partial
u/\partial n=0$ on $S$, which is replaced in Theorem~\ref{thmrelax} by
the weaker assumptions (\ref{div}) and (\ref{divu}). This is due to the
fact that the set $S$ can be approximated by a sequence $S_k$ of sets
with an increasing number of connected components, and this leads to a homogenization process (known as the sieve problem, see,
e.g.,  \cite{Att-Pic}, \cite{Dam}, and \cite{Mur}), where the Neumann condition is replaced in the limit by a transmission condition.
\end{section}

\begin{section}{PROOF OF THE PROPERTIES OF THE SLOPE}

We shall use the following compactness and lower semicontinuity theorem.
\begin{theorem}\label{lsc} Let $E\subset\Om$ be an $\hone$-measurable set with $\hone(E)<+\infty$, and let 
$u_k$ be a sequence in $SBV(\Om)\cap L^\infty(\Om)$ such that $u_k$ is bounded in  $L^\infty(\Om)$, $\hone(S(u_k))$ is bounded, and $\nabla u_k$ is bounded in $L^2(\Om;\R^2)$. 
 Then  there exist a subsequence, still denoted $u_k$, and a function $u\in SBV(\Om)\cap L^\infty(\Om)$, such that 
$u_k\to u$ a.e.\ in $\Om$, $\nabla u_k\wto \nabla u$ weakly in  $L^2(\Om;\R^2)$, and $\hone(S(u)\setmeno E)\le\liminf_k\hone(S(u_k)\setmeno E)$.

If, in addition, there exists a function $\psi\in H^1(\Om)$ such that $u_k-\psi\in SBV^2_0(\Om,E)$ for every $k$, then $u-\psi\in SBV^2_0(\Om,E)$.
\end{theorem}

\begin{proof}
The former statement is proved in \cite{A} (see also \cite{A-F-P}) when $E=\emptyset$. The proof of the general case can be found in \cite[Theorem~2.8]{DM-F-T}. The latter statement can be obtained by considering an extension of all functions to a larger domain.
\end{proof}

We begin by proving Proposition~\ref{necessary}.

\begin{proof}[Proof of Proposition~\ref{necessary}]
Let  $\varphi\in SBV^2_0(\Om,S)$. By the definition of slope (see
(\ref{slopeb})) and of the functional $\F$ (see (\ref{MS})) we have
\begin{equation}\label{phi}
|\partial\F|_b(u,S)\geq
\limsup_{\varepsilon\to0+}\frac{\F(u,S)-\F(u+\varepsilon\varphi,S)}{\varepsilon\|\varphi\|}=-
2\frac{(\nabla u|\nabla\varphi)}{\|\varphi\|}\,.
\end{equation}
Therefore the linear functional $\varphi\mapsto (\nabla u|\nabla\varphi)$
is continuous on $SBV^2_0(\Om,S)$ with respect to the $L^2$ norm. Thus
there exists  $f\in L^2(\Om)$ such that (\ref{nec2}) holds. As
$SBV^2_0(\Om,S)$ is dense in $L^2(\Om)$, (\ref{phi}) implies that
$2\|f\|\leq |\partial\F|_b(u,S)$.

The proof of (\ref{nec1}) is similar.
\end{proof}

We now prove Theorem~\ref{smooth}.
\begin{proof}[Proof of Theorem~\ref{smooth}] We prove the theorem only
for $|\partial\F|_{b}(u,S)$. By Proposition~\ref{necessary} we have
only to prove that $|\partial\F|_b(u,S)\leq 2\|f\|$. Assume, by
contradiction, that
$|\partial\F|_b(u,S)> 2\|f\|$.  Then there exist a constant $\alpha$,
with $0<\alpha\leq+\infty$, and a sequence  $v_k\in SBV(\Om)$ such that
$v_k\to u$ in $L^2(\Om)$, the traces of $v_k$ and $u$ agree on
$\partial\Om$, and
\begin{equation}\label{alpha}
\lim_{k\to\infty}\frac{\|\nabla u\|^2-\|\nabla
v_k\|^2-\hone(S(v_k)\setmeno S)}
{\|u-v_k\|}=\alpha+2\|f\|\,.
\end{equation}
This implies that $\liminf_k(\|\nabla u\|^2-\|\nabla v_k\|^2)\ge0$ and $\limsup_k\hone(S(v_k)\setmeno S)\le0$. Indeed, if this is not the case 
the numerator in (\ref{alpha})
would have a negative limit (along a suitable subsequence), which
contradicts $\alpha>0$.  These inequalities show that $\|\nabla v_k\|$ and 
$\hone(S(v_k)$ are bounded uniformly with respect to~$k$.
By
lower semicontinuity (Theorem~\ref{lsc}) we have also
$\limsup_k(\|\nabla u\|^2-\|\nabla v_k\|^2)\leq0$. Therefore
\begin{equation}\label{h1}
\|\nabla v_k\|\to\|\nabla u\|\qquad\hbox{and}\qquad \hone(S(v_k)\setmeno S)\to0\,.
\end{equation}

As $\Om\setmeno S$ has a finite number of  connected components, the
function $u$ belongs to $L^\infty(\Om)$. By a truncation argument
(changing, if needed, the value of $\alpha$), it is not restrictive to
assume that
\begin{equation}\label{linfty}
\|v_k\|_\infty\leq\|u\|_\infty\,,
\end{equation}
where $\|\cdot\|_\infty$ denotes the $L^\infty$ norm.

Under our hypotheses on $u$ there exists a constant $L$, with $0<L<+\infty$, such that the restriction of $u$ to  each  connected component $U$ of $\Om\setmeno S$ has Lipschitz constant $L$.
Let us fix $\varepsilon>0$, with $2L^2\varepsilon+2L\varepsilon<1$, such that for every $y\in\partial U$ and every $\rho\in{]0,\varepsilon[}$ the set  $\partial B(y,\rho)\cap U$ is connected and $\hone(\partial B(y,\rho)\cap U)\ge 2\rho$.
Let $\Phi_\varepsilon(w):=\|\nabla w\|^2+\frac{1}{\varepsilon}((|f|+\varepsilon)(w-u)|w-u)-2(f|w)$, and let
$W_k$ be the set of all functions $w$ such that
$w-u\in SBV^2_0(\Om,S\cup S(v_k))$.
It is clear that $w_k$ is a solution of the minimum  problem
\begin{equation}\label{wk}
\min_{w\in W_k}\Phi_\varepsilon(w)
\end{equation}
if and only if it minimizes
\begin{equation}\label{wk1}
\min_{w\in W_k}\Big\{\|\nabla w\|^2+\frac{1}{\varepsilon}((|f|+\varepsilon)(w-u-\varepsilon
g_\varepsilon)|w-u-\varepsilon g_\varepsilon)\Big\}\,,
\end{equation}
where $g_\varepsilon:=f/(|f|+\varepsilon)$.
By a truncation argument, we can find a minimizing sequence of
(\ref{wk1}) whose $L^\infty$ norm is bounded by $\|u\|_\infty+\varepsilon$.
Since the approximate gradients in the minimizing sequence are
uniformly bounded in $L^2(\Om;\R^2)$, we can apply the compactness and
lower semicontinuity theorem (Theorem~\ref{lsc}) and we obtain that
the minimum problems (\ref{wk}) and  (\ref{wk1}) have a solution $w_k$,
which is unique by strict convexity.

Let us prove that
\begin{eqnarray}
w_k\to u && \hbox{strongly in }L^2(\Om)\,,\\
\nabla w_k\to \nabla u  && \hbox{strongly in }L^2(\Om;\R^2)\,.
\end{eqnarray}
Indeed,
  $\|w_k\|_\infty$, $\|\nabla w_k\|$, and $\hone(S(w_k))$ are uniformly
bounded, so that by Ambrosio's compactness theorem (Theorem~\ref{lsc})
there exist a subsequence, not relabelled, and a function $ w\in
SBV(\Om)$ such that $w_k\to w$ in $L^2(\Om)$,  $\nabla w_k\wto \nabla
w$ weakly in $L^2(\Om;\Rtwo)$, and $\hone({S(w)\setmeno S})\le
\liminf_k\hone(S(w_k)\setmeno S)$.  Since $S(w_k)\subsethn  S\cup S(v_k)$
and $\hone(S(v_k)\setmeno S)\to0$  we have that
$S(w)\subsethn S$, hence $w\in SBV^2(\Om,S)$.
  As $w_k-u\in SBV^2_0(\Om,S\cup S(v_k))$, we deduce also that
$w-u\in SBV^2_0(\Om,S)$.
Let $v$ be a function such that $v-u\in SBV^2_0(\Om,S)$. As $v\in W_k$,
by the minimality of $w_k$ in  (\ref{wk}) we get
$\Phi_\varepsilon(w_k)\le \Phi_\varepsilon(v)$.
Passing to the limit as $k\to\infty$, the lower
semicontinuity implies that $w$ minimizes $\Phi_\varepsilon$ on the set $W$
of all functions $v$ with $v-u\in SBV^2_0(\Om,S)$. Since,  by
hypothesis, $-\Delta u=f$ in ${\mathcal D}'(\Om)$  and $\partial
u/\partial n=0$ on $S$, the function $u$ minimizes $\|\nabla
v\|^2-2(f|v)$ on $W$. Therefore
$$
\Phi_\varepsilon(u)=\|\nabla u\|^2-2(f|u)\le\|\nabla w\|^2-2(f|w)\le \Phi_\varepsilon(w)
\,,
$$
so that $u=w$ by the uniqueness of the minimizer of $\Phi_\varepsilon$ on $W$.
To prove that the convergence of $\nabla w_k$ to $\nabla u$ is strong,
we observe that $u\in W_k$, so that, by the minimality of $w_k$ in
(\ref{wk}), we have
$\Phi_\varepsilon(w_k)\le \Phi_\varepsilon(u)$. This implies that $\|\nabla
u\|^2\ge\limsup_k\|\nabla w_k\|^2$, which gives the strong convergence.

Since $\Phi_\varepsilon(w_k)\le \Phi_\varepsilon(v_k)$ by (\ref{wk}), we have
\begin{eqnarray}
& \|\nabla u\|^2-\|\nabla v_k\|^2-\hone(S(v_k)\setmeno S)\le
    \|\nabla u\|^2-\|\nabla w_k\|^2-\nonumber\\
  &\displaystyle {}- \frac{1}{\varepsilon}((|f|+\varepsilon)(w_k-u)|w_k-u)+
    \frac{1}{\varepsilon}((|f|+\varepsilon)(v_k-u)|v_k-u)+\nonumber\\
  & {}+ 2(f|w_k-u)- 2(f|v_k-u)-
     \hone(S(v_k)\setmeno S)\,.\nonumber
\end{eqnarray}
We shall prove that for $k$ large enough
\begin{equation}\label{*}
  \|\nabla u\|^2-\|\nabla w_k\|^2-\frac{1}{\varepsilon}((|f|+\varepsilon)(w_k-u)|w_k-u) +
2(f|w_k-u)
\le\hone(S(v_k)\setmeno S)\,,
\end{equation}
hence
\begin{equation}\label{vkwk}
\|\nabla u\|^2-\|\nabla v_k\|^2-\hone(S(v_k)\setmeno S)\le
\frac{1}{\varepsilon}((|f|+\varepsilon)(v_k-u)|v_k-u) - 2(f|v_k-u)\,.
\end{equation}
Since
$\frac{1}{\varepsilon}(|f|+\varepsilon)(v_k-u)\to 0$ strongly in $L^2(\Om)$
as $k\to\infty$ by the dominated convergence theorem, and
$-2(f|v_k-u)\le 2\|f\|\,\|v_k-u\|$,
inequality (\ref{vkwk})  contradicts the fact that $\alpha>0$ in
(\ref{alpha}).

It remains to prove (\ref{*}).
The Euler condition for the minimum problem (\ref{wk}) implies that
\begin{equation}\nonumber
(\nabla w_k|\nabla\varphi)+\frac{1}{\varepsilon}((|f|+\varepsilon)(w_k-u)|\varphi)-(f|\varphi)=0
\end{equation}
for every $\varphi\in SBV^2_0(\Om,S(v_k)\cup S)$. Taking $\varphi=u-w_k$, and using the identity $\|\nabla u\|^2-\|\nabla w_k\|^2=(\nabla u|\nabla u-\nabla  w_k)
+(\nabla w_k|\nabla u-\nabla  w_k)$,
we obtain
\begin{eqnarray*}
& \displaystyle
\|\nabla u\|^2-\|\nabla w_k\|^2-
\frac{1}{\varepsilon}((|f|+\varepsilon)(w_k-u)|w_k-u) +
2(f|w_k-u) =\\
& \displaystyle =(\nabla u|\nabla u-\nabla
w_k) +  (f|w_k-u)\,.
\end{eqnarray*}

Integrating by parts on  each connected component $U$  of
$\Om\setmeno S$ (for a justification under our regularity assumptions see, e.g., the proof of (2.39) in \cite{DM-M-S}) we obtain
$$
(\nabla u|\nabla u-\nabla  w_k) +  (f|w_k-u) =
\int_{S(v_k)\setminus S}\frac{\partial u}{\partial n}[w_k]\,d\hone\,,
$$
so that (\ref{*}) becomes
$$
\int_{S(v_k)\setminus S}\frac{\partial u}{\partial n}[w_k]\,d\hone
\le\hone(S(v_k) \setmeno S)\,.
$$
Since $|\partial u/\partial n|\le L$,
it is enough to show that
\begin{equation}\label{esssup}
\hone\hbox{-}\mathop{\rm ess\,sup}_{S(v_k)\setminus S}|[w_k]|\le \frac1L
\end{equation}
for $k$ large enough.

To do this, we fix a connected component $U$ of $\Om\setmeno S$.
We will prove that for $k$ large enough
there exists a finite number of balls $B(y_i,\rho_i)$, depending possibly on $k$, with
$0<\rho_i<\varepsilon$, such that
\begin{eqnarray}
& \overline U\subset \bigcup_iB(y_i,\rho_i)\label{cover}\,,\\
& \displaystyle \partial B(y_i,\rho_i)\cap
U\cap S(v_k)=\emptyset\,,\label{empty}\\
& \displaystyle \mathop{\hone\hbox{-}\rm ess\,sup}_{U\cap \partial B(y_i,\rho_i)}
|w_k-u|\le \varepsilon\,.\label{esssupi}
\end{eqnarray}
Using an argument related to the maximum principle we will
then show that (\ref{esssupi}) implies that
\begin{equation}\label{osci}
\mathop{\rm osc}_{U\cap  B(y_i,\rho_i)}w_k\le2 L\varepsilon+2\varepsilon\,.
\end{equation}
This gives immediately
\begin{equation}\nonumber
\mathop{\hone\hbox{-}\rm ess\,sup}_{U\cap  B(y_i,\rho_i)\cap S(v_k)}|[w_k]|\le
2L\varepsilon+2\varepsilon\le \frac 1 L\,.
\end{equation}
Since this estimate does not depend on $i$ nor on the connected component $U$,
we obtain~(\ref{esssup}).

We begin by  proving that (\ref{esssupi}) implies (\ref{osci}).
Let $x_i\in U\cap B(y_i,\rho_i)$ and let $m_i:=u(x_i)-L\varepsilon-\varepsilon$ and
$M_i:=u(x_i)+L\varepsilon+\varepsilon$.
Since $u$ is $L$-Lipschitz on $U$ we have
\begin{equation}\label{m0}
m_i+\varepsilon\le u\le M_i-\varepsilon \qquad\hbox{on }U\cap  B(y_i,\rho_i)\,.
\end{equation}
By (\ref{esssupi}) we have
\begin{equation}\label{m2}
m_i\le w_k\le M_i\qquad\hone\hbox{ a.e.\ on }U\cap  \partial B(y_i,\rho_i)\,.
\end{equation}
Moreover, the function $g_\varepsilon$ which appears in (\ref{wk1}) satisfies 
\begin{equation}\label{m1}
m_i\le u+\varepsilon g_\varepsilon\le M_i\qquad\hbox{a.e.\ on }U\cap  B(y_i,\rho_i)\,,
\end{equation}

Let $z_k\in SBV(\Om)$ be defined by
$$
z_k:=\left\{\begin{array}{lr} (m_i\lor w_k)\land M_i &
\hbox{ on }U\cap  B(y_i,\rho_i)\,,\\
\\
w_k & \hbox{ elsewhere}\,.
\end{array}\right.
$$
By (\ref{m2}) we have $z_k\in SBV^2(\Om,S\cup S(v_k))$ and by (\ref{m0})
the traces of $z_k$ and $u$ agree on $\partial\Om$. Therefore $z_k\in W_k$. By
(\ref{m1}) we have
\begin{eqnarray*}
& \displaystyle\|\nabla z_k\|^2+\frac{1}{\varepsilon}((|f|+\varepsilon)(z_k-u-\varepsilon
g_\varepsilon)|z_k-u-\varepsilon g_\varepsilon)\le \\
& \displaystyle \le \|\nabla w_k\|^2+\frac{1}{\varepsilon}((|f|+\varepsilon)(w_k-u-\varepsilon
g_\varepsilon)|w_k-u-\varepsilon g_\varepsilon)\,.
\end{eqnarray*}
Since the solution of (\ref{wk1}) is unique we have $z_k=w_k$, hence
$m_i\le w_k\le M_i$ a.e.\ on $U\cap  B(y_i,\rho_i)$, which implies (\ref{osci}).

It remains to prove that we can find a finite number of balls $B(y_i,\rho_i)$,
with $0<\rho_i<\varepsilon$, satisfying (\ref{cover})--(\ref{esssupi}).
Choose $k$ large enough to
have  
\begin{eqnarray}
& \displaystyle \hone(S(v_k)\setmeno S)<\frac\varepsilon8\,,\label{k1}\\ 
& \displaystyle 64\,(\pi+\frac{1}{\varepsilon^2})\,(\|w_k-u\|^2+ \|\nabla w_k-\nabla u\|^2)<\varepsilon^2\,,\label{k2}
\end{eqnarray}
and define $H_\varepsilon:=\{x\in U:d(x,U^c)\ge  \varepsilon/2 \}$. 
Since $H_\varepsilon$ is compact, it can be covered by a finite number of balls $B(y_i,\varepsilon/4)$ with $y_i\in H_\varepsilon$. 
Note that the balls $B(y_i,\varepsilon/2)$ are contained in $U$.
Let 
$E_i:=\{\rho\in{]\varepsilon/4,\varepsilon/2[}:\partial B(y_i,\rho)\cap S(v_k)=\emptyset\}$.
By (\ref{k1}) we have $\hone(E_i)\ge\varepsilon/8$. Moreover
$$
\int_{E_i}d\rho\int_{\partial B(y_i,\rho)}
|\nabla_\tau w_k-\nabla_\tau u|^2\,d\hone
\le \|\nabla w_k-\nabla u\|^2
$$
and
$$
\int_{E_i}d\rho\int_{\partial B(y_i,\rho)}|w_k- u|^2\,d\hone\le
\|w_k-u\|^2\,,
$$
where $\nabla_\tau$ denotes the tangential gradient.
As $\hone(E_i)\ge\varepsilon/8$,
we may choose $\rho_i\in E_i$ such that
\begin{eqnarray}
& \displaystyle\int_{\partial B(y_i,\rho_i)}|\nabla_\tau w_k-\nabla_\tau u|^2\,d\hone\le \frac{16}{\varepsilon} \|\nabla w_k-\nabla u\|^2\label{rho1}\\
& \displaystyle\int_{\partial B(y_i,\rho_i)}|w_k-u|^2\,d\hone\le
 \frac{16}{\varepsilon}\|w_k-u\|^2\,.\label{rho2}
\end{eqnarray}
Using the one-dimensional estimate
\begin{equation}
\sup_{t\in{]0,a[}}|\zeta(t)|^2\le\frac{2}{a}\int_0^a|\zeta(t)|^2dt+
2a\int_0^a|\zeta'(t)|^2dt\,,\label{1d}
\end{equation}
from  (\ref{k2}), (\ref{rho1}), and (\ref{rho2}) we obtain (\ref{esssupi}).

Let us cover now
$K_\varepsilon:=\{x\in \overline U:d(x,U^c)\le \varepsilon/2\}$. Since
$K_\varepsilon$ is compact we can cover it with a finite number of balls
$B(y_i,3\varepsilon/4)$ with $y_i\in \partial U$. By our choice of $\varepsilon$ 
the set $\partial B(y,\rho)\cap U$ is connected and 
$\hone(\partial B(y_i,\rho)\cap U)\ge \varepsilon$  for $3\varepsilon/4<\rho<\varepsilon$.
 As in the case of balls centered in
interior points of $U$ there exist radii $\rho_i$ with $3\varepsilon/4<\rho_i<\varepsilon$ such that (\ref{empty}), (\ref{rho1}), and (\ref{rho2}) hold with $\partial B(y_i,\rho_i)$ replaced by $\partial B(y_i,\rho_i)\cap U$.  Using again (\ref{1d}) and the lower bound 
$\hone(\partial B(y_i,\rho_i)\cap U)\ge \varepsilon$ we obtain (\ref{esssupi}).
\end{proof}

The following example shows that in the previous theorem it is not enough to assume that $S$ is a manifold with boundary.

\begin{example}\label{example1}
Let $S$ be a closed segment contained in $\Om$  with endpoints $a$, $b$, and let $u\in H^1(\Om\setmeno S)$ be a harmonic function which satisfies the Neumann boundary condition on $S$. It is well known that there exists a constant $\kappa$ such that in a neighbourhood $U$ of $a$ we have
$$
u-\kappa\sqrt{2\rho/\pi}\sin(\theta/2)\in H^2(U\setmeno S)\cap H^{1,\infty}(U\setmeno S) \,,
$$
where $\rho$ and $\theta$ are the polar coordinates around $a$ with $\theta=\pi$ on $S$ (see, e.g., \cite[Theorem 4.4.3.7 and Section~5.2]{Gri1} or \cite[Appendix~1]{MSh}). Moreover, any constant $\kappa$ can be obtained by a suitable choice of the boundary conditions satisfied by $u$ on $\partial\Om$.

For every $s\ge 0$, let $S_s$ be the closed segment obtained by adding to $S$ a collinear segment of length $s$ starting from~$a$, and let
$u(s)$ be the solution of the
problem $\Delta u(s)=0$ on $\Om\setmeno S_s$, $\partial u(s)/\partial n=0$ on $S_s$, and $u(s)=u$ on $\partial\Om$. Then 
\begin{equation}\label{ex10}
|{\partial\F}|_b(u,S)\ge\limsup_{s\to0}\frac{(\F(u(0),S)-\F(u(s),S_s))^+}{\|u(s)-u(0)\|}\,.
\end{equation}
By \cite[Theorem 6.4.1]{Gri2}  the derivative with respect to $s$ of $\F(u(s),S_s)$ at $s=0$ exists and is equal to $1- \kappa^2$. Moreover, a similar argument shows that the function $s\mapsto u(s)$ has a derivative $\dot u(0)$  (in the $L^2$-sense) at $s=0$.
Therefore (\ref{ex10}) gives
$$
|{\partial\F}|_b(u,S)\ge (\kappa^2-1)^+\frac{1}{\|\dot u(0)\|}\,.
$$ 
In particular, if $|\kappa|>1$, we conclude that $|{\partial\F}|_b(u,S)>0$.
\end{example}

The next example shows that $|{\partial\F}|_b$ is not lower semicontinuous.

\begin{example}\label{example2} Under the hypotheses of the previous example, assume in addition that $\Om$ has a $C^2$ boundary. Let $A_k:=\{x\in\Om : d(x,S)<1/k\}$, $S_k=\partial A_k$, and let $\psi_k$ be a sequence of functions in $C^2(\overline\Om)$ which converge to $u$ strongly in $H^1$ near $\partial\Om$. Let  $u_k$ be the solution of $\Delta u_k=0$ in $\Om\setmeno S_k$, which satisfies the Neumann boundary condition $\partial u_k/\partial n=0$ on $S_k$, the Dirichlet boundary condition $u_k=\psi_k$ on $\partial\Om$, and vanishes  on $A_k$. Since $S_k$ is connected and converges to $S$ in the Hausdorff metric (and also in the sense of $\sigma^2$-convergence, see Definition~\ref{sigma2}), we deduce that $u_k$ converges to $u$ strongly  in $L^2(\Om)$ and $\nabla u_k$ converges to $\nabla u$ strongly  in $L^2(\Om;{\bf R}^2)$ (see, e.g., \cite[Theorem~5.1]{DM-T}). By Theorem~\ref{smooth} we have $|\partial\F|_{b}(u_k,S_k)=0$, hence $|\overline{\partial\F}|_b(u,S)=0$, while $|\partial\F|_{b}(u,S)>0$ (see Example~\ref{example1}). This shows that $|\partial\F|_{b}$ is not lower semicontinuous.
\end{example}

\end{section}

\begin{section}{PROOF OF THE PROPERTIES OF THE RELAXED SLOPE}

The definition of relaxed slope is based on the following notion of convergence of sets, introduced in \cite{DM-F-T} and used in the study of quasi-static crack growth. 

\begin{definition}\label{sigma2}
We say that a sequence $S_k$ $\sigma^2$-converges to $S$ if $S_k,\, 
S\subset \Om$, $\hone(S_k)$ is bounded uniformly with respect to $k$,  and the following two conditions are satisfied:
\begin{itemize}
\item[{\rm(a)}] if $u_j\in SBV^2(\Om,S_{k_j})\cap L^\infty(\Om)$ for some sequence $k_j\to\infty$, $u_j$ is bounded in $L^\infty(\Om)$, $u_j\to u$ a.e.\ in $\Om$, and $\nabla u_j\wto\nabla u$ weakly in $L^2(\Om;\R^2)$, then $u\in SBV^2(\Om,S)\cap L^\infty(\Om)$;
\item[{\rm(b)}] there exist a function $u\in SBV^2(\Om,S)\cap L^\infty(\Om)$, with $S(u)\eqhn S$, and a 
sequence $u_k$, bounded in $L^\infty(\Om)$, such that $u_k\in SBV^2(\Om,S_k)$ for every~$k$, $u_k\to u$ a.e.\ in $\Om$, and $\nabla u_k\wto\nabla u$ weakly in $L^2(\Om;\R^2)$.
\end{itemize}
\end{definition}

We begin by proving Proposition~\ref{relaxed}.

\begin{proof}[Proof of Proposition \ref{relaxed}]
Assume that $|\overline{\partial\F}|_b(u,S)<+\infty$.
Then there exist $(u_k,S_k)$ such that $u_k\to u$ strongly in
$L^2(\Om)$, $\nabla u_k\wto \nabla u$ weakly in $L^2(\Om;\Rtwo)$, $S_k$
$\sigma^2$-converges to $S$, and $\lim_k|{\partial\F}|_b(u_k,S_k)=
|\overline{\partial\F}|_b(u,S)$. By Proposition~\ref{necessary}
there exist $f_k\in L^2(\Om)$, with
$2\|f_k\|\le |{\partial\F}|_b(u_k,S_k)$,  such that $(\nabla
u_k|\nabla\varphi)=(f_k|\varphi)$ for every $\varphi\in
SBV^2_0(\Om,S_k)$. It follows that, up to a subsequence, $f_k$ converges
weakly in $L^2(\Om)$ to some function $f$. By lower
semicontinuity
$$
2\|f\|\leq 2\liminf_k\|f_k\|\leq \lim_k|\partial\F|_b(u_k,S_k)=
|\overline{\partial\F}|_b(u,S)\,.
$$
Let $\varphi\in C_0^\infty(\Om)$. As $\varphi\in SBV^2_0(\Om,S_k)$,
we have $(\nabla u_k|\nabla\varphi)=(f_k|\varphi)$, 
and passing to the limit as $k\to\infty$ we get (\ref{div}).
As $\varphi u_k\in SBV^2_0(\Om,S_k)$ we have also
$$
(\nabla u_k|\varphi\nabla u_k)+(\nabla u_k|u_k\nabla\varphi)=(f_k|\varphi u_k)\,.
$$
Passing to the limit as
$k\to\infty$, when $\varphi\ge 0$ we obtain (\ref{divu}) using lower semicontinuity with respect to weak convergence in the first term.

The proof for $|\overline{\partial\F}|(u,S)$ is analogous.
\end{proof}

To prove Theorem \ref{thmrelax} we shall use some properties of the Newtonian capacity and of quasicontinuous representatives of functions in Sobolev spaces, for which we refer to \cite{Eva-Gar} and~\cite{H-K-M}.
Given an orientation of the $C^1$ manifold $S\setmeno F$, for every 
$v\in H^1(\Om\setmeno S)$ the %(pre\-cise values of) 
the traces $v^+$ and $v^-$ on the positive and negative faces of $S\setmeno F$ are 
defined cap-q.e.\ on $S\setmeno F$, hence the jump $[v]:=v^+-v^-$ is defined cap-q.e.\ 
 on $S$ and is cap-quasicontinuous on 
$S\setmeno F$. If $v\in H^1(\Om\setmeno S)\cap L^\infty(\Om\setmeno S)$, by using cut-off functions which vanish in a neighbourhood of $F$, it is easy to prove that
$v\in SBV^2(\Om,S)$ and that 
$Dv\,\lfloor\, S=[v]\,n\,\hone\lfloor\, S$,
%$Dv\,\LL S=[v]\,n\,\hone\LL S$, 
where $n$ is the 
oriented unit normal to $S\setmeno F$. Conversely, if $v\in SBV^2(\Om,S)$ its restriction to $\Om\setmeno S$ belongs to $H^1(\Om\setmeno S)$.
Let 
$$
H^1_{0,\partial\Om}(\Om\setmeno S):=\{v\in H^1(\Om\setmeno S): v=0 \hbox{ on }\partial\Om\}\,.
$$
The previous remarks show that $SBV^2_0(\Om,S)\cap L^\infty(\Om)$ can be identified with $H^1_{0,\partial\Om}(\Om\setmeno S)\cap L^\infty(\Om\setmeno S)$.

Under the assumptions of Theorem \ref{thmrelax}, let $\nu$ be the nonnegative Radon measure on $\Om$ defined by
\begin{equation}\label{nu}
|\nabla u|^2-{\rm div}(u\nabla u)+\nu=f\,u\qquad\hbox{in }{\mathcal
D}'(\Om)\,.
\end{equation}
Since  $|\nabla u|^2$ and $f\,u$ belong to $L^1(\Om)$, while ${\rm div}(u\nabla u)$ belongs to $H^{-1}(\Om)$ it turns out that 
$\nu$ vanishes on all sets of capacity zero.
As $\nabla u $ is the distributional gradient of $u$ on $\Om\setmeno S$, from (\ref{div}) we get $(\nabla u|\nabla (\varphi u))=(f|\varphi u)$ for every $\varphi\in C^\infty_c(\Om\setmeno S)$, which implies that ${\rm supp}\,\nu\subset S$. Using a standard approximation argument we can prove that
\begin{equation}\label{nu2}
(\nabla u|\nabla u\, \varphi) + (u\nabla u|\nabla \varphi)+\int_S\varphi\,d\nu=(f\,u|\varphi) 
\end{equation}
for every $\varphi\in H^1_0(\Om)\cap L^\infty(\Om)$.

Let $\mu$ be the nonnegative Borel measure defined by
\begin{equation}\label{mu}
\mu(B):=\left\{
\begin{array}{lr}
\displaystyle\int_B\frac{d\nu}{[u]^2}  &  \hbox{if }{\rm
cap}(B\cap\{[u]=0\})=0\,,    \\
\\ 
 +\infty & \hbox{otherwise}\,.
\end{array}
\right.
\end{equation}
We note that $\mu$ vanishes on all sets of capacity zero.

\begin{lemma}\label{lemmamu}
Under the assumptions of Theorem~\ref{thmrelax}, let  $\mu$ and $\nu$ be defined by (\ref{mu}) and (\ref{nu}). Then 
$[u]\in L^2(S,\mu)$, $\nu=[u]^2\mu$,  and 
\begin{equation}\label{pmu}
(\nabla u|\nabla v)+([u]|[v])_{S,\mu}=(f|v)
\end{equation}
for every $v\in H^1_{0,\partial\Om}(\Om\setmeno S)$ with 
$[v]\in  L^2(S,\mu)$,
where $(\cdot|\cdot)_{S,\mu}$ denotes the scalar product in
$L^2(S,\mu)$.
\end{lemma}
\begin{proof}
We have
$$
\int_S[u]^2d\mu=\int_{S\cap\{[u]\neq0\}}[u]^2d\mu=
\nu(S\cap\{[u]\neq0\})\le\nu(S)<+\infty\,,
$$
hence $[u]\in L^2(S,\mu)$.

By (\ref{mu}), in order to prove that $\nu=[u]^2\mu$ it is enough to show that 
$\nu(S\cap\{[u]=0\})=0$. Since $\nu(F)=0$ it is enough to show that 
for every $y\in S\setmeno F$ there exists an open neighbourhood $U$ of $y$ such that
\begin{equation}\label{U}
\nu(U\cap S\cap\{[u]=0\})=0\,.
\end{equation}
To this aim we consider a neighbourhood $U$ such that $U\cap S$ is a $C^1$ manifold and $U\setmeno S$ has two connected components $U^\oplus$ and $U^\ominus$, corresponding to the positive and negative faces of $U\cap S$. By possibly reducing $U$ we may assume that there exist two functions $u^\oplus$, $u^\ominus\in H^1(U)\cap L^\infty(U)$ such that $u^\oplus=u$ a.e.\ on  $U^\oplus$
and $u^\ominus=u$ a.e.\ on $U^\ominus$, so that $[u]=u^\oplus-u^\ominus$ cap-q.e.\ on $U\cap S$.

Let $\varphi\in C^\infty_c(U)$ and let $F_\varepsilon:\R\to\R$ be the Lipschitz function defined by
$F_\varepsilon(t):=1+t/\varepsilon$ for $-\varepsilon\le t\le 0$, $F_\varepsilon(t):=1-t/\varepsilon$ for $0\le t\le \varepsilon$, and $F_\varepsilon(t):=0$ for $|t|>\varepsilon$. 
We will prove that 
\begin{equation}\label{intfe}
\int_S F_\varepsilon([u])\,\varphi \,d\nu\to0\,.
\end{equation}
Since $F_\varepsilon([u])\to 1_{\{[u]=0\}}$ pointwise on $S$, by the dominated convergence theorem (\ref{intfe}) implies 
$$
\int_{S\cap\{[u]=0\}}\varphi \,d\nu=0\,,
$$
which gives
(\ref{U}) by the arbitrariness of $\varphi$.

To prove (\ref{intfe}) we use $v:=u^\oplus F_\varepsilon(u^\oplus-u^\ominus)\varphi$ and $w:=F_\varepsilon(u^\oplus-u^\ominus)\varphi$ as test functions in (\ref{div})  and (\ref{nu2}), respectively. We obtain
\begin{equation*}
\begin{array}{c}
\displaystyle(\nabla u| \nabla u^\oplus
F_\varepsilon(u^\oplus-u^\ominus)\varphi)+ 
(u^\oplus\nabla u| F'_\varepsilon(u^\oplus-u^\ominus)(\nabla u^\oplus-\nabla u^\ominus) 
\varphi)+
{}\\
\displaystyle{} \vphantom{\int_S} +(u^\oplus\nabla u| \nabla \varphi
F_\varepsilon(u^\oplus-u^\ominus)) =(f|u^\oplus F_\varepsilon(u^\oplus-u^\ominus)\varphi)
\end{array}
\end{equation*}
and 
\begin{equation*}
\begin{array}{c}
\displaystyle (\nabla u| \nabla u
F_\varepsilon(u^\oplus-u^\ominus)\varphi)+ 
(u\nabla u| F'_\varepsilon(u^\oplus-u^\ominus)(\nabla u^\oplus-\nabla u^\ominus) 
\varphi)+
{}\\
\displaystyle{} +(u\nabla u| \nabla \varphi
F_\varepsilon(u^\oplus-u^\ominus)) +\int_S F_\varepsilon([u])\varphi d\nu =(f|uF_\varepsilon(u^\oplus-u^\ominus)\varphi)\,.
\end{array}
\end{equation*}
Subtracting term by term we get
\begin{equation}\label{diff}
\begin{array}{c}
\displaystyle(\nabla u|(\nabla u^\oplus-\nabla u)
F_\varepsilon(u^\oplus-u^\ominus)\varphi)+ 
((u^\oplus-u)\nabla u| F'_\varepsilon(u^\oplus-u^\ominus)(\nabla u^\oplus-\nabla u^\ominus) 
\varphi)+
{}\\
\displaystyle{} +((u^\oplus-u)\nabla u| \nabla \varphi
F_\varepsilon(u^\oplus-u^\ominus)) -\int_S F_\varepsilon([u])\varphi d\nu=(f|(u^\oplus-u) F_\varepsilon(u^\oplus-u^\ominus)\varphi)\,.
\end{array}
\end{equation} 
Since $u^\oplus=u$ and $\nabla u^\oplus=\nabla u$ a.e.\ on the set $\{u^\oplus=u^\ominus\}$, the first, the third term and the last term in (\ref{diff}) tend to $0$ by the dominated convergence theorem. As $|(u^\oplus-u)F'_\varepsilon(u^\oplus-u^\ominus)|<1$ and the measure of the set $\{0<|u^\oplus-u^\ominus|<\varepsilon\}$ tends to $0$, the second term tends to $0$ too. This proves (\ref{intfe}), which concludes the proof of the fact that $\nu=[u]^2\mu$.

To prove (\ref{pmu}) we need the following density result.
\end{proof}

\begin{lemma}\label{density}
Under the assumptions of Theorem~\ref{thmrelax} for every 
$v\in H^1_{0,\partial\Om}(\Om\setmeno S)$ with 
$[v]\in  L^2(S,\mu)$, there exist $\varphi_k$, $\psi_k\in C^\infty_c(\Om)$ such that, setting
$v_k:=\varphi_k+\psi_ku$, we have
\begin{eqnarray}
& v_k\to v\qquad\hbox{strongly in }H^1(\Om\setmeno S)\,,\label{(1)}\\
& [v_k]\to[v]\qquad\hbox{strongly in }L^2(S,\mu)\,.\label{(2)}
\end{eqnarray}
\end{lemma}

\begin{proof}[Proof of Lemma \ref{lemmamu} (continuation)] 
It remains to prove (\ref{pmu}). Let $v\in H^1_{0,\partial\Om}(\Om\setmeno S)$ with 
$[v]\in  L^2(S,\mu)$ and let $\varphi_k$, $\psi_k$ be as in Lemma~\ref{density}. Taking 
$\varphi_k$ and $\psi_k$ as test function in (\ref{div}) and (\ref{nu2}), respectively, and using the equality $\nu=[u]^2\mu$, we obtain 
\begin{eqnarray*}
& (\nabla u|\nabla\varphi_k)=(f|\varphi_k)\\
& (\nabla u|\psi_k\nabla u)+(\nabla u|u \nabla \psi_k)+([u]|\psi_k[u])_{S,\mu}=(f|u\psi_k)\,.
\end{eqnarray*}
Adding term by term we get 
$$
(\nabla u|\nabla v_k)+([u]|[v_k])_{S,\mu}=(f|v_k)\,,
$$
where $v_k:=\varphi_k+\psi_ku$. Passing to the limit thanks to Lemma~\ref{density} we obtain (\ref{pmu}).
\end{proof}

\begin{proof}[Proof of Lemma~\ref{density}]  
Since every function $v\in H^1_{0,\partial\Om}(\Om\setmeno S)$ with 
$[v]\in  L^2(S,\mu)$ can be approximated by truncations, it is enough to prove the lemma when $v$ is bounded. Since $F$ has capacity zero, we may also assume that $v$ vanishes a.e.\ in a neighbourhood of $F$. 
By using a partition of unity, it is enough to prove the lemma in one of the following cases:
\begin{itemize}
\item[(a)] ${\rm supp}\, v\cap S=\emptyset$;
\item[(b)] ${\rm supp}\, v$ is contained in an open set $U$ such that $U\cap S$ is a $C^1$ manifold and $U\setmeno S$ has two connected components $U^\oplus$ and $U^\ominus$, corresponding to the positive and negative faces of $U\cap S$.
\end{itemize}
In the former case it is enough to take $\psi_k=0$, and the result follows from the density of $C^\infty_c(\Om\setmeno S)$ in $H^1_0(\Om\setmeno S)$.

In the latter case there exist two functions $v^\oplus$, $v^\ominus\in H^1_0(U)\cap L^\infty(U)$, with compact support in $U$, such that $v^\oplus=v$ a.e.\ on  $U^\oplus$
and $v^\ominus=v$ a.e.\ on $U^\ominus$, so that $[v]=v^\oplus-v^\ominus$ cap-q.e.\ on $U\cap S$.

Since $[v]\in L^2(S,\mu)$, from the definition of $\mu$ (see (\ref{mu})) we obtain 
${\rm cap}(\{{[v]\neq 0}\}\cap\{{[u]=0}\})=0$.  
It is not restrictive to assume that there exists $\varepsilon>0$ such that $\{{[v]\neq 0}\}\subset \{{|[u]|>\varepsilon}\}$ (see Lemma~\ref{vke} below).

Let $u^\oplus$ and $u^\ominus$ be the functions defined in the proof of Lemma~\ref{lemmamu}, and let $G_\varepsilon:\R\to\R$ be the Lipschitz function defined by
$G_\varepsilon(t):=1/t$ for $|t|\ge \varepsilon$ and $G_\varepsilon(t):=t/\varepsilon^2$ for $|t|\le \varepsilon$.
Since $G_\varepsilon(u^\oplus-u^\ominus)\in H^1(U)\cap L^\infty(U)$, the function 
$(v^\oplus-v^\ominus)\,G_\varepsilon(u^\oplus-u^\ominus)$ belongs to $H^1_0(U)\cap L^\infty(U)$ and coincides with $[v]\,G_\varepsilon([u])$ cap-q.e.\ on $S$. Since $\{[v]\neq 0\}\subset \{|[u]|>\varepsilon\}$,  we have $[u]\,[v]\,G_\varepsilon([u])=[v]$
cap-q.e.\ on $S$. 

By the density of  $C^\infty_c(U)$ in $H^1_0(U)$ there exists a sequence  $\psi_k\in C^\infty_c(U)$ which is  bounded in $L^\infty(U)$ and converges to $(v^\oplus-v^\ominus)\,G_\varepsilon(u^\oplus-u^\ominus)$ strongly in $H^1_0(U)$. In particular, up to a subsequence, $\psi_k\to [v]\,G_\varepsilon([u])$ cap-q.e.\ on $S$. Hence 
$\psi_k[u]\to [v]$ strongly in $L^2(S,\mu)$. 
Let $w_k\in H^1(U\setmeno S)$ be defined by
$w_k:=(u^\oplus-u^\ominus)\,(v^\oplus-v^\ominus)\,G_\varepsilon(u^\oplus-u^\ominus)-\psi_k(u^\oplus-u^\ominus)$
 on $U^\oplus$, and $w_k:=0$ on $U^\ominus$. 
 Note that $w_k\to 0$ strongly in $H^1(U\setmeno S)$ and ${\rm supp}\,w_k\subset\subset U$.
 As $[w_k]=[v]-\psi_k[u]$  cap-q.e.\ on $S$, 
the function $v-\psi_k u-w_k$ has no jump on $S$, hence it belongs to $H^1(U)$. Since ${\rm supp}(v-\psi_k u-w_k)\subset\subset U$, we have $v-\psi_k u-w_k\in H^1_0(U)$, thus there exists a sequence $\varphi_k\in C^\infty_c(U)$ such that $v-\psi_k u-w_k-\varphi_k\to 0$ strongly in $H^1_0(U)$. It is then clear that the sequence $v_k$ defined in the statement of the lemma satisfies (\ref{(1)}) and (\ref{(2)}).
\end{proof}

\begin{lemma}\label{vke}
Let $u$ be as in Theorem~\ref{thmrelax}, and  let $U$, $U^\oplus$, $U^\ominus$ be as in Lemma~\ref{density}. Let $v\in H^1(U\setmeno S)\cap L^\infty(U)$ with  ${\rm supp}\,v\subset\!\subset U$ and $[v]\in L^2(S,\mu)$.
Then there exist a sequence $v_k\in H^1(U\setmeno S)\cap L^\infty(U)$ and a sequence $\varepsilon_k>0$ such that 
\begin{eqnarray}
&& v_k\to v\quad\hbox{strongly in }H^1(U\setmeno S)\,,\label{v1}\\
&& [v_k]\to[v]\quad\hbox{strongly in }L^2(S\cap U,\mu)\,,\label{v2}\\
&& \{[v_k]\neq 0\}\subset\{|[u]|>\varepsilon_k\}\,.\label{v3}
\end{eqnarray}
\end{lemma}

\begin{proof}
As we noticed in the proof of Lemma~\ref{density}, the definition of $\mu$ (see (\ref{mu})) implies that ${\rm cap}(\{[v]\neq 0\}\cap\{[u]=0\})=0$.  

Let us prove that there exists a sequence $w_k\in H^1(U\setmeno S)\cap L^\infty(U)$ such that 
\begin{eqnarray}
&& w_k\to v\quad\hbox{strongly in }H^1(U\setmeno S)\,,\label{111}\\
&& [w_k]\to[v]\quad\hbox{strongly in }L^2(S\cap U,\mu)\,,\label{112}\\
&& |[w_k]|\le k|[u]|\quad\hbox{cap-q.e.\ on }S\,.\label{113}
\end{eqnarray}
Let $g_k\colon U\to \R$ be the function defined by
$g_k:=(v^\oplus-v^\ominus)\lor 0$ cap-q.e.\ on $U\setmeno S$ and $g_k:=[v]\lor 0\land k|[u]|$ cap-q.e.\ on $S$. The functions $g_k$ are cap-quasi lower semicontinuous, the sequence $g_k$ is increasing and converges to (the cap-quasi continuous representative of) $(v^\oplus-v^\ominus)\lor 0$ cap-q.e.\ on $U$. 
By \cite[Lemma 1.6]{DM} there exists a sequence $w_k^{(1)}\in H^1_0(U)$ such that $0\le w_k^{(1)}\le |[v]|\land k|[u]|$ cap-q.e.\ on $S\cap U$ and $w_k^{(1)}\to (v^\oplus-v^\ominus)\lor 0$ strongly in $H^1_0(U)$.
Similarly, there exists a sequence  $w_k^{(2)}\in H^1_0(U)$ such that $-( |[v]|\land k|[u]|)\le w_k^{(2)}\le 0$ and $w_k^{(2)}\to (v^\oplus-v^\ominus)\land 0$ strongly in $H^1_0(U)$. Then $w_k^{(1)}+w_k^{(2)}\to  v^\oplus-v^\ominus$ strongly in $H^1_0(U)$, hence, up to a subsequence, $w_k^{(1)}+w_k^{(2)}\to [v]$ cap-q.e.\ on $S\cap U$. Moreover $|w_k^{(1)}+w_k^{(2)}|\le |[v]|\land k|[u]|$ cap-q.e.\ on $S\cap U$.
Let us define $w_k^\oplus:=v^\ominus+w_k^{(1)}+w_k^{(2)}$ and $w_k^\ominus:=v^\ominus$. It is then clear that the function $w_k$ defined by
$w_k:=w_k^\oplus$ in $U^\oplus$ and 
$w_k:=w_k^\ominus$ in $U^\ominus$
satisfies (\ref{111})--(\ref{113}).

For every $\varepsilon>0$ let  $T_\varepsilon:\R\to\R$ be the Lipschitz function defined by $T_\varepsilon(t):=t-\varepsilon$, for $t\ge\varepsilon$, $T_\varepsilon(t):=t+\varepsilon$ for $t\le-\varepsilon$, and $T_\varepsilon(t):=0$ for $|t|\le\varepsilon$. Let $w_{k,\varepsilon}\in H^1_{0,\partial U}(U\setmeno S)$ be the function defined by
$$
w_{k,\varepsilon}:=\left\{\begin{array}{lr}
w_k^\ominus+T_{k\varepsilon}(w_k^\oplus-w_k^\ominus) & \hbox{in }U^\oplus\,,\\
\\
w_k^\ominus & \hbox{in }U^\ominus\,.
\end{array}\right.
$$
Then $[w_{k,\varepsilon}]=T_{k\varepsilon}([w_k])$ cap-q.e.\ on $S\cap U$, $w_{k,\varepsilon}\to w_k$ strongly in $H^1(U\setmeno S)$, and $[w_{k,\varepsilon}]\to [w_k]$ strongly in $L^2(S\cap U,\mu)$.

We choose now $\varepsilon_k>0$ such that $\|w_{k,\varepsilon_k}-w_k\|+\|\nabla w_{k,\varepsilon_k}-\nabla w_k\|<1/k$ and $\|[w_{k,\varepsilon_k}]-[w_k]\|_{\mu}<1/k$. Then the function $v_k:=w_{k,\varepsilon_k}$ satisfies (\ref{v1}) and (\ref{v2}). Since $[v_k]=T_{k\varepsilon_k}([w_k])$ and $|[w_k]|\le k|[u]|$ cap-q.e.\ on $S\cap U$, we obtain that (\ref{v3}) is satisfied.
\end{proof}

Let $w\in H^1(\Om)$. For every compact subset $M$ of $\Om$ we consider the functional $\G^w_M:L^2(\Om)\to[0,+\infty]$ defined by
\begin{equation}\label{gm}
\G^w_M(v):=\left\{\begin{array}{lr}
\|\nabla v\|^2+\|v\|^2 & \hbox{ if }v-w\in H^1_{0,\partial\Om}(\Om\setmeno M)\,,\\
\\
+\infty & \hbox{otherwise.}
\end{array}\right.
\end{equation}
For every nonnegative Borel measure $\mu$ on $\Om$, vanishing on all sets of capacity zero and  with ${\rm supp}\,\mu\subset S$, we consider the functional $\F_\mu:L^2(\Om)\to[0,+\infty]$ defined by
\begin{equation}\label{fmu}
\F^w_\mu(v):=\left\{\begin{array}{lr}
\|\nabla v\|^2+\|v\|^2+\|[v]\|^2_{S,\mu} & \hbox{ if }v-w\in H^1_{0,\partial\Om}(\Om\setmeno S)\,,\\
\\
+\infty & \hbox{otherwise.}
\end{array}\right.
\end{equation}

For the definition and properties of $\Gamma$-convergence we refer to \cite{DMb} and \cite{Bra}.

\begin{lemma}\label{gamma}
Let $S$ be as in Theorem~\ref{thmrelax}, and let $\mu$ be a nonnegative Borel measure on $\Om$ vanishing on all sets of capacity zero and with ${\rm supp}\,\mu\subset S$. Then there exists a sequence $M_k$ of $C^1$ manifolds with boundary  such that $M_k\subset S$ and 
$\G^{w_k}_{M_k}$ $\Gamma$-converges to $\F^w_\mu$ in $L^2(\Om)$ for every sequence $w_k$ converging to $w$ strongly in $H^1(\Om)$.
\end{lemma}

\begin{proof} This proof is obtained by adapting \cite[Theorem 4.16]{DM-Mos} and \cite[Theorem 2.38]{Att}. First of all we observe that it is enough to prove the lemma when $w_k=w=0$. The general case can be obtained by modifying the functions near the boundary in order to match the boundary conditions.

Let $\mathcal X$ be the set of lower semicontinuous functionals $\G\colon L^2(\Om)\to[0,+\infty]$ with $\G\ge \G^0_S$. Under our regularity assumptions on $S$ we can apply Rellich's Theorem and we obtain that  the set $\{v\in L^2(\Om): \G^0_S(v)\le t\}$  is compact in $L^2(\Om)$ for every $t<+\infty$. Therefore 
$\Gamma$-convergence in $L^2(\Om)$ is metrizable on $\mathcal X$
by \cite[Theorem 10.22]{DMb}. 

Let $\mathcal Y$ be the set of functionals $\G^0_M$, where $M$ runs over all $C^1$ manifolds with boundary contained in $S$. Since $\mathcal Y\subset\mathcal X$, to prove the lemma for $w_k=w=0$ we have to show that $\F^0_\mu\in \overline{\mathcal Y}$, where 
$\overline{\mathcal Y}$ denotes the closure of $\mathcal Y$ with respect to the metric which induces  $\Gamma$-convergence on $\mathcal X$. Indeed, in this case there exists a sequence $M_k$ of $C^1$ manifolds with boundary  such that $M_k\subset S$  and $\G^0_{M_k}$ $\Gamma$-converges to $\F^0_\mu$ in $L^2(\Om)$.

Using the techniques developed in \cite{Cor} it is possible to prove that $\F^0_\mu\in \overline{\mathcal Y}$ whenever $\mu=\hone\lfloor E$ and $E\subset S\setmeno F$ is an arc, i.e., a connected $C^1$ manifold with or without boundary. More in general, the same techniques show that  $\F^0_\mu\in \overline{\mathcal Y}$ whenever $\mu=g\hone\lfloor S$ and $g$ is a step-function, i.e., $g=\sum c_i1_{E_i}$ where $c_i\ge 0$ and $E_i$ are disjoint arcs contained in $S\setmeno F$. 

To continue the proof we need the following lemma.
\end{proof}

\begin{lemma}\label{convmu}
Let $S$ be as in Theorem~\ref{thmrelax}, and
let $\mu_k$, $\mu$ be positive measures in $H^{-1}(\Om)$ with support in $S$. 
If $\mu_k \to\mu$ strongly in $H^{-1}(\Om)$ then $\F^0_{\mu_k}$ $\Gamma$-converges to $\F^0_\mu$ in $L^2(\Om)$.
\end{lemma}

\begin{proof}
Let $v_k$ be a sequence of functions such that $v_k\to v$ strongly in $L^2(\Om)$ and 
$\F^0_{\mu_k}(v_k)\le C$ for some constant $C<+\infty$. Then $v_k$ converges to $v$ weakly in $H^1_{0,\partial\Om}(\Om\setmeno S)$. For every $\varepsilon>0$ let $S_\varepsilon$ be the difference between $S$ and the union of the closed balls with centers in $F$ and radius~$\varepsilon$. As $S_\varepsilon$ is a $C^1$ manifold, arguing as in the proof of Lemma~\ref{density}, for every $\varepsilon$ we can construct $z_k$, $z\in H^1_0(\Om)$ (possibly depending on $\varepsilon$) such that
$z_k=[v_k]$ cap-q.e.\/ on $S_\varepsilon$, $z=[v]$ cap-q.e.\/ on $S_\varepsilon$, and $z_k\wto z$ weakly in $H^1_0(\Om)$. For every $t\in\R$ the sequence $(|z_k|\land t)^2$ converges to $(|z|\land t)^2$ weakly in $H^1_0(\Om)$.  Therefore
$$
\int_{S_\varepsilon} (|z|\land t)^2d\mu\le\liminf_{k\to\infty}\int_{S_\varepsilon} (|z_k|\land t)^2d\mu_k\le
\liminf_{k\to\infty}\int_{S_\varepsilon} z_k^2d\mu_k\,.
$$
Passing to the limit first as $t\to\infty$ and then as $\varepsilon\to0$ we obtain 
$$
\int_{S} [v]^2d\mu\le \liminf_{k\to\infty}\int_{S} [v_k]^2d\mu_k\,,
$$
which implies  
$$
\F^0_\mu(v)\le\liminf_{k\to\infty} \F^0_{\mu_k}(v_k)\,.
$$

Arguing as in the beginning of the proof of Lemma~\ref{density} we see that it is enough to prove the $\Gamma$-limsup inequality when 
$v\in H^1_{0,\partial\Om}(\Om\setmeno S)$ is bounded and vanishes a.e.\/ in a neighbourhood of $F$. In this case there exists $z\in H^1_0(\Om)\cap L^\infty(\Om)$ such that $z=[v]$ cap-q.e.\/ on $S$.  As $z^2\in H^1_0(\Om)$ we have 
$$
\|[v]\|^2_{S,\mu}=\int_Sz^2d\mu=\lim_{k\to\infty}\int_Sz^2d\mu_k=\lim_{k\to\infty}\|[v]\|^2_{S,\mu_k}\,,
$$ 
which implies
$$
\F^0_\mu(v)=\lim_{k\to\infty} \F^0_{\mu_k}(v)\,.
$$
This concludes the proof of the $\Gamma$-limsup inequality.
\end{proof}

\begin{proof}[Proof of Lemma~\ref{gamma} (continuation)]
Since every positive  measure  in $H^{-1}(\Om)$ with support in $S$ can be approximated strongly in $H^{-1}(\Om)$ by measures of the form $\mu_k=g_k\hone\lfloor S$, with $g_k$ step-functions, from the previous step of the proof and from Lemma~\ref{convmu} we deduce that $\F^0_\mu\in \overline{\mathcal Y}$ for every  positive measure $\mu$ in $H^{-1}(\Om)$ with support in $S$. 

Let $\mu$ be a nonnegative Borel measure on $\Om$ vanishing on all sets of capacity zero and with ${\rm supp}\,\mu\subset S$.
By \cite[Lemma 4.15]{DM-Mos} there exist a positive measure $\mu_0$  in 
$H^{-1}(\Om)$ with support in $S$, and a Borel function $g:S\to[0,+\infty]$ such that
$$
\int_S v^2d\mu=\int_S v^2 g\,d\mu_0
$$ 
for every $v\in H^1(\Om)$.

By localizing the problem to an open set $U$ satisfying condition (b) in the proof of Lemma~\ref{density}, and using the functions $v^\oplus$, $v^\ominus$ introduced in that proof we obtain that $\int_S [v]^2d\mu=\int_S [v]^2 g\,d\mu_0$ for every $v\in H^1(\Om\setmeno S)$. Therefore $\F^0_\mu=\F^0_{g\mu_0}$. Let $g_k:=g\lor k$. Then
$g_k\mu_0\in H^{-1}(\Om)$, hence $\F^0_{g_k\mu_0}\in  \overline{\mathcal Y}$. As $\F^0_{g_k\mu_0}$ is increasing and converges pointwise to $\F^0_{g\mu_0}=\F^0_\mu$ we conclude that $\F^0_{g_k\mu_0}$ $\Gamma$-converges to $\F^0_\mu$ in $L^2(\Om)$, hence  $\F^0_\mu\in  \overline{\mathcal Y}$.
\end{proof}

\begin{proof}[Proof of Theorem \ref{thmrelax}]
By Proposition~\ref{relaxed} it is enough to prove that
$|\overline{\partial\F}|_b(u,S)\le 2\|f\|$. To do this we show that,
under the regularity hypotheses of the theorem, there exist $(u_k,S_k)$
with $S_k$ smooth, such that $S_k$ $\sigma^2$-converges to $S$,  
$u_k\to u$ strongly in $L^2(\Om)$, $\nabla u_k\wto \nabla u$ weakly in
$L^2(\Om;\Rtwo)$,  $|{\partial\F}|_b(u_k,S_k)=2\|f_k\|$, and $f_k\to f$ strongly
in $L^2(\Om)$.

Let $w$ be a function in $H^1(\Om)$ having the same trace as $u$ on $\partial\Om$, and let $w_k$ be a sequence in $C^\infty(\overline\Om)$ such that $w_k\to w$ strongly in  $H^1(\Om)$.
Using Lemma~\ref{gamma} we obtain that
 there exists a sequence $M_k$ of $C^1$ manifolds with boundary  such that $M_k\subset S$ and $\G^{w_k}_{M_k}$ $\Gamma$-converges to $\F^w_\mu$ in $L^2(\Om)$, where $\mu$ is the measure given by (\ref{mu}).
  
As (\ref{pmu}) is the Euler equation for the problem
\begin{equation}\label{minu}
\min_{v\in L^2(\Om)} \Big\{\F^w_\mu(v)-2(f+u|v)\Big\}\,, 
\end{equation}
the function $u$ is the solution of (\ref{minu}).
Let $v_k$ be the minimizer of
\begin{equation}\label{mink}
\min_{v\in L^2(\Om)} \Big\{\G^{w_k}_{M_k}(v)-2(f+u|v)\Big\}\,. 
\end{equation} 
By $\Gamma$-convergence $v_k\to u$ strongly in $L^2(\Om)$. As $\nabla v_k$ are bounded in $L^2(\Om;\Rtwo)$, we have that $\nabla v_k\wto\nabla u$ weakly in $L^2(\Om;\Rtwo)$.

We now approximate each manifold $M_k$ in the Hausdorff metric by a sequence $M_k^\varepsilon$, $\varepsilon>0$, of $C^\infty$ manifolds without boundary having the same number of connected components as $M_k$, and we consider the  
 solutions $v_k^\varepsilon$ of the minimum problems
\begin{equation}\label{minek}
\min_{v\in L^2(\Om)}\Big\{\G^{w_k}_{M_k^\varepsilon}(v)-2(f+u|v)\Big\}\,. 
\end{equation} 
By \cite[Theorems 2.3 and 4.2]{DM-Eb-Pon} $v_k^\varepsilon\to v_k$ strongly in $L^2(\Om)$ and 
$\nabla v_k^\varepsilon\to \nabla v_k$ strongly in $L^2(\Om;\Rtwo)$.
Therefore we can choose $\varepsilon_k>0$ such that setting $u_k:=v_k^{\varepsilon_k}$ and $S_k:=M_k^{\varepsilon_k}$ we have $\|u_k-v_k\|<1/k$, $\|\nabla u_k-\nabla v_k\|<1/k$, and the Hausdorff distance between $S_k$ and $M_k$ is less then $1/k$. 

Then $u_k\to u$ strongly in $L^2(\Om)$,  $\nabla u_k\wto\nabla u$ weakly in $L^2(\Om;\Rtwo)$. 
The Euler equation of (\ref{minek}) implies that 
\begin{equation}\label{eulerk}
\left\{\begin{array}{lr}
-\Delta u_k=f+u-u_k & \hbox{on }\Om\setmeno S_k\,,\\
\frac{\partial u_k}{\partial n}=0 & \hbox{on }S_k\,,\\
u_k={w_k} & \hbox{on }\partial\Om\,.
\end{array}\right.
\end{equation} 
The regularity theory for (\ref{eulerk}) (see, e.g., \cite[Theorem~3.17]{Tr}) gives that $u_k\in W^{2,p}(\Om\setmeno S_k)$. In particular, for every connected component $U$ of $\Om\setmeno S_k$ we have that $u_k\in C^1(\overline U)$. By Theorem~\ref{smooth} we have $|\partial\F|_b(u_k,S_k)=2\|f+u-u_k\|$. 

Therefore it remains to prove that $S_k$ $\sigma^2$-converges to $S$. The condition on the Hausdorff distance between $S_k$ and $M_k$, together with the inclusion $M_k\subset S$, implies that $S_k$ is contained in the closed 
$\varepsilon$-neighbourhood $S^{(\varepsilon)}$ of $S$ for $k$ large enough.
If $z_k\in SBV^2(\Om)\cap L^\infty(\Om)$ is a sequence as in  condition (a) in Definition~\ref{sigma2} for $S_k$, then $z_k$ is bounded in $H^1(\Om\setmeno S^{(\varepsilon)})$. This implies that its limit $z$ belongs to  $H^1(\Om\setmeno S^{(\varepsilon)})$, hence $S(z)\subsethn S^{(\varepsilon)}$. As $\varepsilon$ is arbitrary, we deduce that $S(z)\subsethn S$, so that condition (a) is satisfied. As for condition (b), it is enough to take $u$ and a suitable truncation of $u_k$. 
\end{proof}

\end{section}

\bigskip

\noindent {\bf Acknowledgments.} { This work is part of the Project
``Calculus of Variations" 2002, supported by the Italian Ministry of
Education, University, and Research.}
\bigskip
\bigskip

{
}

\end{document}